\documentclass[graybox]{svmult}

\usepackage{mathptmx}
\usepackage{helvet}
\usepackage{courier}
\usepackage{type1cm}
\usepackage{makeidx}
\usepackage{graphicx}
\usepackage{multicol}
\usepackage[bottom]{footmisc}

\usepackage{amscd}
\usepackage{amsfonts}
\usepackage{amssymb}
\usepackage{latexsym}
\usepackage[arrow, matrix, curve]{xy}

\newcommand{\lkadncm}{\newcommand}

\newtheorem{lkad_theorem}{Theorem}[section]
\newtheorem{lkad_prop}[lkad_theorem]{Proposition}
\newtheorem{lkad_lemma}[lkad_theorem]{Lemma}
\newtheorem{lkad_cor}[lkad_theorem]{Corollary}

\def\lkadN{\mathbb{N}\,}
\def\lkadK{\mathbb{K}\,}

\lkadncm{\lkadEnd}{\mbox{\rm End}\,}
\def\lkadHom{\mbox{\rm Hom}\,}
\def\lkadrad{\mbox{\rm rad}\,}
\def\lkad|{\, | \,}
\def\lkadh{\, \sim \,}

\def\lkado{\otimes}
\def\lkadbra{\langle}
\def\lkadket{\rangle}
\def\lkadIndec{\mbox{\rm Indec}\,}
\def\lkadeps{\varepsilon}

\begin{document}

\title*{Subalgebra depths within the path algebra of an acyclic quiver}
\author{Lars Kadison and Christopher J. Young}
\institute{Lars Kadison \at  Departamento de Matematica, Faculdade de Ciencias da Universidade do Porto, 687 Rua Campo Alegre, 4169-007 Porto, Portugal, \email{lkadison@fc.up.pt}
\and Christopher J. Young \at Departamento de Matematica, Faculdade de Ciencias da Universidade do Porto, 687 Rua Campo Alegre, 4169-007 Porto, Portugal, \email{booloon\_b457@hotmail.com}}


\maketitle

\abstract*{Constraints are given on the depth of diagonal subalgebras in generalized
triangular matrix algebras.  The depth of the top subalgebra $B \cong A \lkadrad A$ in
a finite, connected, acyclic quiver algebra $A$ over an algebraically closed field $\lkadK$ is then computed.  Also the depth of the primary arrow subalgebra  $1\lkadK + \lkadrad A = B$ in $A$ is obtained. The two types
of subalgebras have depths $3$ and $4$ respectively, independent of the number of vertices.  An upper bound on depth is obtained for the quotient of a subalgebra pair.}

\abstract{Constraints are given on the depth of diagonal subalgebras in generalized
triangular matrix algebras.  The depth of the top subalgebra $B \cong A / \lkadrad A$ in
a finite, connected, acyclic quiver algebra $A$ over an algebraically closed field $\lkadK$ is then computed.  Also the depth of the primary arrow subalgebra  $1\lkadK + \lkadrad A = B$ in $A$ is obtained. The two types
of subalgebras have depths $3$ and $4$ respectively, independent of the number of vertices.  An upper bound on depth is obtained for the quotient of a subalgebra pair.}

\section{Introduction}

Given a subalgebra pair, one extracts a (minimum) depth from a comparison of $n$-fold tensor products
of the subalgebra pair with one another in a meaningful way.  The interesting case is when an $(n+1)$-fold tensor product divides a multiple of the $n$-fold tensor product in the sense of Krull-Schmidt unique factorization into indecomposable bimodules, or more generally as a bimodule isomorphism with a direct summand.
The bimodule structures on the $n$-fold tensor products are naturally any one of four possibilities
as left and right modules over the subalgebra or overalgebra.  The least restrictive of these conditions
is two-sided over the subalgebra and we fix the depth in the situation mentioned above to be
$2n+1$; for mixed bimodules, we have the left and right depth $2n$ conditions \cite{BDK}.  The most
stringent condition, as bimodules of the overalgebra, is H-depth $2n-1$ \cite{LKH}, and is useful to ordinary depth gauging as well
when the overalgebra has nice bimodules such as a separable algebra (see Proposition~\ref{lkadison/prop-A1A} below).

Comparing the tensor-square of an algebra extension with the overalgebra as mixed bimodules leads to a characterization of the Galois extension \cite{KS, K2008, LKJPAA}.  Thus not unexpectedly
the depth two condition placed on Hopf subalgebras is equivalent to the normality condition
with respect to the adjoint actions \cite{BK}.   The depth three condition is satisfied by a subalgebra $B \subseteq A$ when, in a suitably nice category of bimodules, $A$ contains all $B^e$-indecomposables  that can possibly appear up to isomorphism in decompositions of tensor products $A \otimes_B \cdots \otimes_B A$ \cite{BK, LK}.  Semisimple complex subalgebra pairs
of each depth $n \in \lkadN$ are noted in \cite{BKK} via bipartite graphs and inclusion matrices for $K_0(B) \rightarrow K_0(A)$.

In the paper \cite{BDK} it was shown that the depth of a finite group algebra extension is bounded
by twice the index of the normalizer of the subgroup in the group.
In the papers \cite{BuK, BKK, BDK, D, F} the depth of certain group algebra extensions are computed; for example, \cite{F} computes the depth of all the subgroups of $PSL(2,q)$ viewed as complex group algebras.  In \cite{BKK} the complex group algebras associated to the permutation groups are shown
to have depth $d(S_n, S_{n+1}) = 2n-1$; in \cite{BDK}, this same result is shown to not depend on the ground ring.

It was noted in the paper \cite{LK} that a subalgebra $B$ in a finite-dimensional algebra
$A$ has finite depth $d(B,A)$ if $B^e$ has finite representation type; below we note that this holds
 if $A^e$ has finite representation type.  In addition it is possible in algebras without involution that a subalgebra having left depth $2n$ may not have right depth $2n$.
Moreover, the matrix power inequality characterizing depth $n$ subalgebra pairs of semisimple complex algebras in \cite{BuK, BKK} breaks down in the presence of indecomposables of length greater than one.  For these reasons, it becomes interesting to begin a study of depth of
subalgebras in path algebras of quivers.   A reasonable place to start is with acyclic quivers for whose path algebras there is a classic theorem about which have finite representation type in terms of Dynkin diagrams and the underlying graphs \cite{ASS}.
This paper computes the depth of the top and arrow subalgebras of the path algebra of a finite, connected, acyclic quiver.  In Section~3 we note constraints on the depth of a diagonal subalgebra of a generalized matrix ring. We also note an inequality of depth in case the subalgebra contains ideals of the overalgebra, perhaps useful in computing depth of certain  subalgebras of bounded quiver algebras.  In the last Section~6 of concluding remarks we discuss other subalgebras of certain quiver algebras and their depth.

\section{Preliminaries on depth}

Given a unital associative ring $R$ and unital $R$-modules $M$ and $N$, we say that $M$ divides
$N$ and write $M \lkad| N$ if $N \cong M \oplus *$ as $R$-module for some (unnamed) complementary
module. If there are natural numbers $r$ and $s$ such that $N \lkad| rM = M \oplus \cdots \oplus M$ and $M \lkad| sN$, then $M$ and $N$ are  H-equivalent (or similar), as $R$-modules; denoted by $M \lkadh N$. Note that this is  indeed an equivalence relation.  In this case their endomorphism rings $\lkadEnd M_R$ and $\lkadEnd N_R$
are Morita equivalent with Morita context bimodules $\lkadHom (M_R,N_R)$ and
$\lkadHom (N_R,M_R)$ (with module actions and Morita pairings given by composition).

  If $M$ and $N$ are in a category of finitely generated $R$-modules having unique factorization into
indecomposables, then $M$ and $N$ have the same indecomposable constituents if and only if
$M$ and $N$ are H-equivalent modules.  If $F$ is an additive endofunctor of the category of $R$-modules, then $M \lkadh N$ implies $F(M) \lkadh F(N)$; which in practice means that H-equivalent bimodules
may replace one another in certain H-equivalences of  tensor products. In addition, $M \lkadh N$
and $U \lkadh V$ implies $M \oplus U \lkadh N \oplus V$.

 Throughout this paper, let $A$ be a unital associative ring and $B \subseteq A$ a subring where $1_B =  1_A$.  Note the natural bimodules ${}_BA_B$ obtained by restriction of the natural $A$-$A$-bimodule (briefly $A$-bimodule) $A$, also  to the natural bimodules ${}_BA_A$, ${}_AA_B$ or ${}_BA_B$, which
are referred to with no further notation.  Equivalently we denote the proper ring extension $A \supseteq B$ occasionally by $A \lkad| B$.  (Often results are valid as well for a ring homomorphism $B \rightarrow A$
and its induced bimodules on $A$.)

Let $C_0(A,B) = B$, and for $n \geq 1$,
$$ C_n(A,B) = A \lkado_B \cdots \lkado_B A \ \ \ \mbox{\rm ($n$ times $A$)} $$
For $n \geq 1$, the $C_n(A,B)$ has a natural $A$-bimodule  structure given by
$a(a_1 \otimes \cdots \otimes a_n)a'  = aa_1 \otimes \cdots \otimes a_na'$.
Of course, this bimodule structure restricts to $B$-$A$-, $A$-$B$- and $B$-bimodule structures as we may need them.
Let $C_0(A,B)$ denote the natural $B$-bimodule $B$ itself.
Recall from \cite{BDK, LK} that a subring $B \subseteq A$ has  right depth $2n$
if \begin{equation}
\label{lkadison/eq: subringdepth}
C_{n+1}(A,B) \lkadh C_n(A,B)
\end{equation}
as natural $A$-$B$-bimodules; left depth $2n$ if the same condition holds as $B$-$A$-bimodules;
if both left and right conditions hold, it has depth $2n$;
and depth $2n+1$ if the same condition holds as $B$-bimodules.  If condition~(\ref{lkadison/eq: subringdepth}) holds in its strongest form as $A$-$A$-modules for $n \geq 1$ the subring $B \subseteq A$ is said to
have H-depth $2n-1$; H-depth is investigated in \cite{LKH}.

Note that if the subring
has left or right depth $2n$, it automatically has depth $2n+1$ by restriction to $B$-bimodules.
Also note that if the subring has depth $2n+1$, it has depth $2n+2$ by
tensoring the H-equivalence by $- \lkado_B A$ or $A \lkado_B -$. The \textit{minimum depth} (or just depth when the context makes it clear) is denoted
by $d(B,A)$; if $B \subseteq A$ has no finite depth, write $d(B,A) = \infty$.  There is hidden in this a
subtlety: if there is a subring $B \subseteq A$ of left depth $2n$ but not of right depth
$2n$, then it has depth $2n+1$, left and right depth $2n+2$, and nevertheless its minimum depth is $2n$.  There is not a published example of such a subring at present (but a search for this
must occur outside the class of QF extensions \cite[Th.\ 2.4]{LK}).  Note too that
if $B \subseteq A$ has H-depth $2n-1$, it has depth $2n$ by restriction.

In practice one only need check half of the condition in (\ref{lkadison/eq: subringdepth}) to establish depth
$2n$ or $2n+1$ of a ring extension $A \supseteq B$.  This is due to the fact that it is always the
case that $C_n(A,B) \lkad| C_{n+1}(A,B)$ for $n \geq 1$ via appropriate face and degeneracy maps
in the relative homological bar complex; e.g. the $A$-$A$-epimorphism $a_1 \otimes a_2 \mapsto a_1a_2$ is split by the $B$-$A$-monomorphism
$a \mapsto 1 \otimes_B a$, whence $C_1(A,B) \lkad| C_2(A,B)$ as $B$-$A$-bimodules.

For a $k$-algebra $B$ let $B^e$ denote $B \otimes_k B^{\rm op}$. For a finite
dimensional algebra $A$ let $n_A$ denote the cardinal number of isomorphism classes of indecomposable finitely generated $A$-modules.  Of course each of the $B^e$-modules $C_n(A,B)$ are finitely generated when $A$ is a finite dimensional algebra.
\begin{lkad_prop}
\label{lkadison/prop-A1A}
Let $B \subseteq A$ be a subring pair of finite dimensional algebras.  If  $B^e$ has finite representation type, then $d(B,A) \leq 1 + 2n_{B^e} $. If $A^e$ has finite representation type, then $d(B,A) \leq 2n_{A^e}$.  If $A \otimes B^{\rm op}$ has finite representation type, then $d(B,A) \leq 2n_{A \otimes B^{\rm op}}$.
\end{lkad_prop}
\begin{proof}
If $B^e$ has finite representation type, it is shown in \cite{LK} that subring depth $d(B,A)$ is finite based on two basic facts.  First, a finitely generated module $M$ over a finite dimensional algebra
divides a multiple of another module $N$ if and only if their Krull-Schmidt unique factorization into indecomposable modules possess
the  indecomposable constituents satisfying $\lkadIndec (M) \subseteq \lkadIndec (N)$;
then $M$ and $N$ are H-equivalent iff $\lkadIndec (M) = \lkadIndec (N)$.  Secondly, from  $C_n(A,B) \lkad| C_{n+1}(A,B)$ we obtain $\lkadIndec C_n(A,B)$ as sequence of subsets of a finite number of indecomposables that grows with $n$.

If $A^e$ has finite representation type, then one applies the same argument with growing $\lkadIndec C_n(A,B)$,
this time as $A$-$A$-bimodules, which shows that $C_{N+1}(A,B)$ and $C_N(A,B)$ are H-equivalent
after at most $N = n_{A^e}$ steps.  Then the minimum H-depth $d_H(B,A) \leq 2N-1$, and one
notes by restricting modules that $d(B,A) \leq 2N$.   The last statement is proven similarly using the definition of even depth.
\end{proof}

\begin{lkad_cor}
Suppose $B \subseteq A$ is a subalgebra pair where either $A$ or $B$ is a separable algebra.
Then depth $d(B,A)$ is finite.
\end{lkad_cor}


\section{Constraints on subring depth in triangular matrix rings}
\label{lkadison/sect-triangular}

Let $R$ and $S$ be unital associative rings.  Suppose ${}_SM_R$ is a unital $S$-$R$-bimodule
as suggested by the notation.  There is a triangular matrix ring, denoted by $A$, associated with this data,
\begin{equation}
A := \left( \begin{array}{cc}
R & 0 \\
M & S
\end{array}
\right)
\end{equation}
with the obvious matrix addition and multiplication,
which defines a well-known class of examples in the demonstration of independence of axioms in ring
theory such as left and right noetherian property of rings.

Note the subring of diagonal matrices in $A$
is isomorphic (and identified) with $R \times S$.  The obvious split epimorphism of rings $ A \rightarrow R \times S$ is denoted by $\pi: \left( \begin{array}{cc} r & 0 \\
m & s  \end{array} \right) \mapsto (r, s)$.  The mapping $\pi$ is of course an isomorphism if $M = 0$.    Also note the orthogonal idempotents $e_1 = (1_R,0)$ and $e_2 = (0,1_S)$, where $A = e_1A \oplus e_2 A e_1 \oplus  A e_2$.

Let $R'$ be a unital subring of $R$, and $S'$ a unital subring of $S$.  Then
$B := R' \times S' $ is a subalgebra of diagonal matrices in $A$.  We will be interested
in the depth $d(B,A)$.  At first we will dispose of the case $M = 0$ and note that
$d(R' \times S', R \times S) = \max \{ d(R',R), d(S',S) \}$.  (This proposition should be
compared with \cite[Prop.\ 3.15]{BKK}.)

\begin{lkad_prop}
The depth of a subalgebra of a direct product of rings is given by
$$ d(R' \times S', R \times S) = \max \{ d(R',R), d(S',S) \}. $$
\end{lkad_prop}
\begin{proof}
Let $A = R \times S$ and $B = R' \times S'$.
 Note that the central orthogonal idempotents $e_1, e_2 \in B \subseteq A$.  It follows that there is the following isomorphism of $n$-fold
tensor products (any $n \in \lkadN$),
\begin{equation}
\label{lkadison/eq: decomp}
C_n(A,B) \cong C_n(R,R') \oplus C_n(S,S')
\end{equation}
 as $B$-$B$-, $A$-$B$- and $B$-$A$-bimodules up to a trivial extension of
for example $R$-module to $A$-module by $S \cdot x = 0$, all elements $x$ in the module. Such a
decomposition holds as well for bimodule homomorphisms between $n$- and $n+1$-fold tensor products.

Let $2m+1 \geq \max \{ d(R',R), d(S',S) \}$.  Then the righthand-side of (\ref{lkadison/eq: decomp}) where $n = m+1$ divides a multiple of the $m$-fold tensor
product of the same form, then so does the lefthand-side.  Hence $d(B,A) \leq 2m+1$. If both depths $d(R',R)$ and $d(S',S)$ are even,
the same argument replacing $2m+1$ with $2m$ suffices to establish $d(B,A) \leq \max \{ d(R',R), d(S',S) \}$. Note that the argument works for $0$-fold tensor product and depth one case too. The reverse inequality follows from applying the central idempotents to $C_n(A,B) \lkadh C_{n+1}(A,B)$.
\end{proof}

 Next  we continue the notation $B = R' \times S'$ and $A$
as the triangular matrix ring formed from the rings $R$, $S$ and the bimodule ${}_SM_R \neq 0$.
Let $\mathcal{M}$ denote a category of modules or bimodules, where left and right subscripts denote the rings in action.
\begin{lkad_lemma}
As abelian categories, $${}_B\mathcal{M}_B \cong {}_{R'}\mathcal{M}_{R'} \oplus {}_{R'}\mathcal{M}_{S'}
\oplus {}_{S'}\mathcal{M}_{R'} \oplus {}_{S'}\mathcal{M}_{S'}$$
\end{lkad_lemma}
\begin{proof}
This isomorphism is induced on objects by ${}_BV_B \mapsto e_1Ve_1 \oplus e_1Ve_2 \oplus e_2Ve_1 \oplus e_2Ve_2$.  Conversely, an object $(W_1, W_2,W_3, W_4)$ on the right side is sent to a matrix
$\left( \begin{array}{cc}
W_1 & W_2 \\
W_3 & W_4
\end{array}
\right)$ with left action by row vectors $(r,s)$ and right action by column vectors $\left( \begin{array}{c}
r' \\
s'
\end{array}
\right)$.  A $B$-bimodule homomorphism $f: V \rightarrow W$ commutes with $e_1,e_2$
from left and right, so that $f$ sends $e_iVe_j$ into $e_iWe_j$ for all $i,j = 1,2$.  Conversely,
a morphism of $2 \times 2$ matrices as before commutes with row and column vectors, and so is a $B$-bimodule homomorphism.
\end{proof}
We now apply the lemma to the $B$-bimodules, the $n$-fold tensor products of the triangular matrix ring $A$ over the diagonal subalgebra $B$.
\begin{lkad_lemma}
For integer $n \geq 1$, $e_1C_n(A,B)e_1 = C_n(R,R')$, $e_1 C_n(A,B) e_2 = 0$, and \\ $e_2C_n(A,B) e_2 = C_n(S, S')$; also
\begin{equation}
\label{lkadison/eq: etwo}
e_2 C_n(A,B) e_1 =  \sum_{r = 0}^{n-1} \oplus \ C_r(S,S') \otimes_{S'} M \otimes_{R'} C_{n-1-r}(R,R')
 \end{equation}
\end{lkad_lemma}
\begin{proof}
For $a_1,\ldots, a_n \in A$, the computations follow from $e_1a_1 \otimes_B \cdots \otimes_B a_n$ $=$ $e_1a_1 e_1 \otimes \cdots \otimes_B a_n$ $= \cdots$ $= e_1a_1 \otimes_B \cdots \otimes_B  e_1a_n$;
moreover, $a_1 \otimes_B \cdots \otimes_B a_n e_2$ $=$ $a_1 \otimes_B \cdots \otimes_B e_2 a_n e_2$
$= \cdots $ $= a_1 e_2 \otimes_B \cdots \otimes_B a_ne_2$; furthermore,
$e_1a_1 \otimes_B \cdots \otimes_B a_n e_2 = 0$ by referring to the last computation and noting
$e_1Ae_2 = 0$.  Naturally, $C_n(e_1A,B) = C_n(R,R')$ since $B = R' \times S'$ and $S'$ acts as zero, so the relative tensor product is given by factoring out by only the nonzero relations; the same is true of $C_n(Ae_2,B) = C_n(S,S')$.

Finally, the last equation follows from $e_2a_1 \otimes_B \cdots \otimes_B a_n e_1$ $=$ $(e_2a_1e_2 + e_2a_1 e_1) \otimes_B \cdots \otimes_B (e_2a_ne_1 + e_1a_n e_1)$ $=\cdots$ $=$
$ \sum_{i=1 }^n a_1 e_2 \otimes_B
\cdots \otimes e_2 a_i e_1 \otimes_B \cdots \otimes_B e_1 a_n$.  This follows from cancellations of the type $\cdots \otimes a_i e_1 \otimes_B \cdots \otimes_B e_2a_j \otimes_B \cdots = 0$ since
$e_1a_k = e_1a_ke_1$, $a_k e_2 = e_2 a_k e_2$ for all $a_k \in A$ and of course $e_1e_2 = 0$.
\end{proof}
Let $d_{\mbox{\rm odd}}(B,A)$ be the smallest odd number greater than or equal to $d(B,A)$,
which we call the odd depth of the subring $B \subseteq A$.
If the depth is finite and already odd, then $d_{\mbox{\rm odd}}(B,A) = d(B,A)$, and otherwise $d_{\mbox{\rm odd}}(B,A) =  d(B,A) + 1$.  In other words, a ring extension $A \lkad| B$
has $d_{\mbox{\rm odd}}(B,A) = 2n+1$ if the natural $B$-$B$-bimodules $C_{n+1}(A,B) \lkadh C_n(A,B)$
and $n$ is the smallest such natural number.
\begin{lkad_theorem}
The odd depth $d_{\mbox{\rm odd}}(B,A)$ satisfies the inequalities,
\begin{equation}
d(B, R \oplus S)  \leq d_{\mbox{\rm odd}}(B,A) \leq d_{\mbox{\rm odd}}(R',R) + d_{\mbox{\rm odd}}(S',S) + 1
\end{equation}
\end{lkad_theorem}
\begin{proof}
If $B \subseteq A$ has depth $2n+1$, then there is $q \in \lkadN$ such that $C_{n+1}(A,B) \oplus V \cong qC_n(A,B)$ for some $B$-$B$-bimodule $V$.  It follows that $e_iC_{n+1}(A,B)e_i \oplus e_iVe_i \cong qe_iC_n(A,B)e_i$ for $i = 1,2$, so that $C_{n+1}(R,R') \lkad| qC_n(R,R')$ and
$C_{n+1}(S,S') \lkad| $  $ q C_n(S,S')$.  It follows that $R' \subseteq R$ and $S' \subseteq S$ both
have depth $2n+1$.  Then $\max \{ d(R',R),$ $ d(S',S) \} \leq d_{\mbox{\rm odd}}(B,A)$.
This completes the proof of the first of the two inequalities.

Next let $R' \subseteq R$ and $S' \subseteq S$ have depths $2n+1$ and $2m + 1$ respectively.
This means that for each integer $s \geq 1$ and $r \geq 0$ there is $q \in \lkadN$ such that $C_{n+s}(R,R') \lkad| qC_{n+r}(R,R')$ as $B$-$B$-bimodules (and similarly for $S' \subseteq S$).
Consider $C_{n+m+2}(A,B)$ as a natural $B$-$B$-bimodule.  By the lemma, $C_{n+m+2}(A,B) \cong$
$$ C_{n+m+2}(R,R') \oplus C_{n+m+2}(S,S') \oplus \sum_{i=0}^{n+m+1}
\oplus C_i(S,S') \otimes_{S'} M \otimes_{R'} C_{n+m+1-i}(R,R') $$
which divides as $B$-$B$-bimodules (due to the depth hypotheses) a multiple of
$$ C_{n+m+1}(R,R') \oplus C_{n+m+1}(S,S') \oplus \sum_{j=0}^{n+m} C_j(S,S') \otimes_{S'} M \otimes_{R'} C_{n+m - j}(R,R'), $$ which is  isomorphic to a multiple of $ C_{n+m+1}(A,B)$. Hence $B \subseteq A$ has depth $2(n+m+1)+1 = 2n+ 2m+3$.  This establishes
that $d(B,A) \leq d_{\mbox{\rm odd}}(B,A) \leq d_{\mbox{\rm odd}}(R',R) + d_{\mbox{\rm odd}}(S',S) + 1$.
\end{proof}
Note that the proof shows that if $R' \subseteq R$ and $S' \subseteq S$ are subrings of finite
depth, then so is $B \subseteq A$, and conversely.

\subsection{Quotient algebras and depth bounds}
Let $B \subseteq A$ be an arbitrary algebra extension and let $I \subseteq B$ be an $A$-ideal. For purposes of expedient notation we write $B_I := B/I$ and similarly for $A_I$. The main purpose of this section is to give some depth bounds for $B_I \subseteq A_I$ as another algebra extension. It turns out that if $d(B,A)$ is finite, then so is $d(B_I, A_I)$.

Recall that if the extension $B \subseteq A$ has odd depth $2n+1$ (even depth $2n$) then $$C_{n+1}(A,B) \sim C_n(A,B)$$ as $B$-bimodules ($A$-$B$- and $B$-$A$-bimodules), which is in general equivalent to saying that there're two $B$-$B$-homomorphisms ($A$-$B$- and $B$-$A$-homomorphisms) $f: C_{n+1}(A,B) \rightarrow m C_n(A,B)$ and $g: m C_n(A,B) \rightarrow C_{n+1}(A,B)$ such that $g \circ f = id$.

\begin{lkad_lemma}[$\pi$ and $\sigma$ properties]\label{cyoung/piprops}
Suppose that $B \subseteq A$ and $I \subseteq B$ are as above. We define the following maps:
\begin{eqnarray*}
\pi &:& C_{n}(A,B) \rightarrow C_{n}(A_I , B_I)      \\
   &:& a_1\otimes \ldots \otimes a_{n} \mapsto \overline{a_1} \otimes \ldots \otimes \overline{a_{n}}. \\ \\
\sigma &:& C_{n+1}(A,B) \rightarrow C_{n+1}(A_I,B_I) \\ &:& a_1\otimes \ldots \otimes a_{n+1} \mapsto \overline{a_1} \otimes \ldots \otimes \overline{a_{n+1}}
\end{eqnarray*}

\noindent These two maps are well-defined and will be $k$-linear as well as satisfying $$\pi(x\heartsuit y) = \overline{x}\pi(\heartsuit)\overline{y} \mbox{ and } \sigma(x\diamondsuit y)= \overline{x}\sigma(\diamondsuit)\overline{y},$$ $\forall x,y \in A$, $\forall \heartsuit \in C_{n}(A,B)$ and $\forall \diamondsuit \in C_{n+1}(A,B)$.
\end{lkad_lemma}

\noindent As will be necessary in our next result we "raise $\pi$ to the $m^{\small \mbox{th}}$ power" in that we define $ \pi' : mC_{n}(A,B) \rightarrow mC_{n}(A_I,B_I)$ in the obvious way: $$(\heartsuit_i) \mapsto (\pi(\heartsuit_i)).$$  The important thing to note however is that $\pi' (x \heartsuit_i y) = \overline{x}\pi'(\heartsuit_i)\overline{y}$, where $x,y \in A$ and $\heartsuit_i \in mC_n(A,B)$, furthermore $\pi'$ is $k$-linear over elements of $mC_n(A,B)$.

\begin{lkad_theorem}
Suppose that $B\subseteq A$ is an algebra extension with depth $2n+1$ ($2n$), suppose also that $I \subseteq B \subseteq A$ is an $A$-ideal. Then $B_I \subseteq A_I$ also has depth $2n+1$ ($2n$). Indeed we can say $d(B_I,A_I) \leq d(B,A)$.

\begin{proof}
\begin{rm}
We prove the odd case because it involves $B$-bimodules and the proof can be extended to the even case with $A$-$B$-bimodules. First, because $B \subseteq A$ has depth $2n+1$ we have $B$-bimodule maps $f: C_{n+1}(A,B) \rightarrow mC_{n}(A,B)$ and $g:  mC_{n}(A,B) \rightarrow C_{n+1}(A,B)$ such that $g\circ f = id$, where $m\geq 1$. We'd like first to find a $B_I$-bimodule map $$\widetilde{f} : C_{n+1}(A_I,B_I) \rightarrow mC_{n}(A_I,B_I) $$ and secondly another $B_I$-bimodule map $$\widetilde{g}: mC_n(A_I, B_I) \rightarrow C_{n+1}(A_I, B_I)$$ such that $\widetilde{g}\circ\widetilde{f} = id$. This enforcing the depth $2n + 1$ condition on $B_I \subseteq A_I$.

We define $\widetilde{f}$ as follows:
\begin{equation}\label{cyoung/fnew}
\widetilde{f}(\overline{a_1} \otimes \ldots \otimes \overline{a_n}) := \pi'\circ f(a_1 \otimes \ldots \otimes a_n)
\end{equation}

\noindent We must show that $\widetilde{f}$ is well-defined, and to that end with some $1 \leq p \leq n$ let $\overline{a_p} = \overline{y}$, that is $a_p = y + t,$ for $t \in I$. Thus
\begin{eqnarray*}
\widetilde{f}(\overline{a_1} \otimes \ldots \otimes \overline{a_p} \otimes \ldots \otimes \overline{a_n}) &=& \pi' f(a_1 \otimes \ldots \otimes y + t \otimes \ldots \otimes a_n) \\ &=& \pi' f(a_1 \otimes \ldots \otimes y \otimes \ldots \otimes a_n) + \pi' f(a_1 \otimes \ldots \otimes t \otimes \ldots \otimes a_n) \\ &=&
\pi' f((a_1 \otimes \ldots \otimes y \otimes \ldots \otimes a_n)) \\ &=& \widetilde{f}(\overline{a_1} \otimes \ldots \otimes \overline{y} \otimes \ldots \otimes \overline{a_n})
\end{eqnarray*}
since $\pi' f(a_1 \otimes \ldots \otimes a_{p-1} \otimes t \otimes a_{p+1} \otimes \ldots \otimes a_n) = \pi' f(a_1 \otimes \ldots \otimes t_1 \otimes 1 \otimes a_{p+1} \otimes \ldots \otimes a_n)$ etc until we have $\pi' (t_p f(1\otimes \ldots \otimes 1 \otimes a_{p+1} \otimes \ldots \otimes x_n)) = \overline{t_p} (\pi' f (1\otimes \ldots\otimes  a_n)) = 0$ (where each $t_i \in I$). This all follows because $I \subseteq B$ is an $A$-ideal with the properties of lemma (\ref{cyoung/piprops}) in effect. Repeating such a process over all $1 \leq p \leq n$  the map will be well-defined.

Now we describe $\widetilde{g}$:

\begin{equation}\label{cyoung/gnew}
\widetilde{g}((\overline{a_1} \otimes \ldots \otimes \overline{a_{n+1}})_i) := \sigma\circ g((a_1 \otimes \ldots \otimes a_{n+1})_i)
\end{equation}

\noindent Proving that $\widetilde{g}$ is well-defined is so similar to the (\ref{cyoung/fnew}) case it can be considered a minor exercise. Furthermore we should notice that $\widetilde{g} \circ \pi' = \sigma \circ g$ straight off. Using (\ref{cyoung/fnew}) and (\ref{cyoung/gnew}) we demonstrate that $\widetilde{g} \circ \widetilde{f} = id$:
\begin{eqnarray*}
\widetilde{g} \circ \widetilde{f}(\overline{a_1} \otimes \ldots \otimes \overline{a_n}) &=& \widetilde{g}\circ \pi'\circ f(a_1 \otimes \ldots \otimes a_n) \\ &
=& \sigma \circ g \circ f(a_1 \otimes \ldots \otimes a_n) \\ &=& \sigma \circ id(a_1 \otimes \ldots \otimes a_n)\\ &=& \overline{a_1} \otimes \ldots \otimes \overline{a_n}
\end{eqnarray*}
\end{rm}
\end{proof}

\end{lkad_theorem}

\begin{lkad_cor}
Given a chain of $A$-ideals $J_0 \subseteq J_1 \subseteq \ldots \subseteq B$ we have $$1 \leq \ldots \leq d(B_{J_1}, A_{J_1}) \leq d(B_{J_0}, A_{J_0}) \leq d(B,A)$$

\begin{proof}
\begin{rm}
The second isomorphism theorem tells us that $(B/ J_0)/(J_1 / J_0) \cong B/J_1$. Apply our last theorem to see that the depth of $(B/J_0)/(J_1 / J_0) \subseteq (A/J_0)/(J_1 / J_0)$ is less than or equal to the depth of $(B/J_0) \subseteq (A/J_0)$, but then we're done.
\end{rm}
\end{proof}
\end{lkad_cor}

\section{Depth of top subalgebra in path algebra of acyclic quiver}
\label{lkadison/sec: top}
Let $Q = (V,E,s,t)$ denote a finite connected acyclic quiver with vertices $V$ of cardinality $|V| = n$
and oriented edges $E$ such that $|E| < \infty$, where an oriented edge or arrow is denoted by $\alpha: a \rightarrow b$, or $(a|\alpha |b) \in E$, where $a = s(\alpha)$ and $b = t(\alpha)$ define the source and target mappings $E \rightarrow V$, respectively.  Since $Q$ is acyclic, there is no loop in $E$, i.e., no arrow $\beta \in E$ such that
$s(\beta) = t(\beta)$; moreover, there are no other cycles, i.e.,  paths $( a | \alpha_1,\ldots,\alpha_r | a )$
of length $r >1$ beginning at a vertex $a$ and ending there (where all $\alpha_i \in E$ and
$s(\alpha_{i+1}) = t(\alpha_i)$, $i = 1,\ldots,r-1$).

Let $\lkadK$ be an algebraically closed field and  let $A = \lkadK Q$ be the path algebra on the quiver $A$  \cite{ASS, P}  with basis the set of all paths, including stationary paths denoted
by $\lkadeps_a = (a||a)$ for each $a \in V$, such that the product of two basis elements is given
by the following concatenation formula:
\begin{equation}
(a| \alpha_1,\ldots,\alpha_r|b)(c|\beta_1,\ldots,\beta_s|d) = \delta_{bc}(a|\alpha_1,\ldots,\alpha_r,\beta_1,\ldots,\beta_s|d).
\end{equation}
The product on $A$ is given by this formula and linearization, which clearly makes $A$ into a graded algebra where $A_s$ denotes the $\lkadK$-vector subspace spanned by paths of length $s$,  a complete set of  primitive orthogonal idempotents are $\{ \lkadeps_a | a \in V \} \in A_0$ and
the radical ideal is $\lkadrad A = A_1 \oplus A_2 \oplus \cdots$, also known as the arrow ideal.

There is always a numbering of the vertices from $1,\ldots,n$ such that $(i | \alpha | j) \in E$
implies $i > j$ \cite[cor.\ 8.6]{P}.  The vertex $n$ is then a source and $1$ a sink.
With such a numbering the algebra $A = \lkadK Q$ is embeddable in a lower triangular matrix algebra
\cite[Lemma 1.12]{ASS}
of the form,
\begin{equation}
\label{lkadison/eq: kaycue}
A = \left( \begin{array}{cccc}
 \lkadeps_1 (\lkadK Q)\lkadeps_1 & 0 & \cdots & 0 \\
\lkadeps_2(\lkadK Q)\lkadeps_1 & \lkadeps_2 (\lkadK Q) \lkadeps_2 & \cdots & 0 \\
\vdots &  \vdots  &  & \vdots \\
\lkadeps_n (\lkadK Q) \lkadeps_1 & \lkadeps_n (\lkadK Q) \lkadeps_2 & \cdots & \lkadeps_n (\lkadK Q) \lkadeps_n
\end{array}
\right)
\end{equation}

Note that $\lkadeps_i (\lkadK Q)\lkadeps_i \cong K$ for each $i = 1,\ldots,n$ since there are no cycles. For example,
if the quiver $Q$ has no multiple arrows between vertices and its underlying graph is a tree,
then there is at most one path between two points $i > j$, so that $\dim \lkadeps_i(\lkadK Q)\lkadeps_j \leq 1$,
and $A = \lkadK Q$ is isomorphic to a subalgebra of the full triangular matrix algebra $T_n(\lkadK) =
\sum_{n \geq i \geq j \geq 1} \lkadK e_{ij}$ (in terms of matrix units $e_{ij}$).

Another example: if $Q = (V,E)$ where $V = \{ 1,2 \}$ and $E = \{ \alpha, \beta: 2 \rightarrow 1 \}$,
then \begin{equation}
\label{lkadison/eq: kronecker-path-alg}
A = \lkadK Q = \left( \begin{array}{cc}
\lkadK & 0 \\
\lkadK^2 & \lkadK
\end{array}
\right)
\end{equation}
From the result of the previous section, we note that with $M = \lkadK^2$, and $B = \lkadK \lkadeps_1 + \lkadK \lkadeps_2$, the depth of $B$ in $A$ is bounded by
\begin{equation}
1 \leq d(B,A) \leq 3.
\end{equation}

For this algebra, one constructs from nilpotent Jordan blocks of order $m$ an infinite sequence of  indecomposable
$A$-modules \cite[pp. 75-76]{ASS}, a tame Kronecker algebra \cite[V111.7]{ARS}.   The algebra $A = \lkadK Q$ has finite representation type
if and only if the underlying (multi-) graph of $Q$ is one of the Dynkin diagrams
$A_n (n \geq 1), D_n (n\geq 4), E_6, E_7, E_8$:
see for example \cite[Gabriel's Theorem, 5.10]{ASS} or \cite[VIII.5.2]{ARS}.

Coming back to the algebra $A$ in (\ref{lkadison/eq: kaycue}), note that $A$ has $n$ augmentations
$\rho_i: A \rightarrow \lkadK$ given by $\rho_i(\lambda_1,\ldots,\lambda_n) = \lambda_i$.  Let $A^+_i$ denote
$\ker \rho_i$, and for a subalgebra $B \subseteq A$, let $B^+_i$ denote $\ker \rho_i \cap B$.
Denote the $n$ $A$-simples of dimension one by ${}_{\rho_i}\lkadK $, and the $n^2$ $A^e$-simples
by $\lkadK_{ij}$ where $a \cdot 1 \cdot b = \rho_i(a) \rho_j(b) 1$ for all $a,b \in A$ and $i,j=1,\ldots,n$.
We have the following
\begin{lkad_lemma}
 Suppose $B \subseteq A$ is a subalgebra of an algebra with augmentations $\rho_1,\ldots,\rho_n$. If $B \subseteq A$ has
right depth $2$, then $AB_i^+ \subseteq B_i^+A$ for each $i = 1,\ldots,n$.  If $B \subseteq A$ has
left depth $2$, then $B_i^+ A \subseteq AB_i^+$ for each $i = 1,\ldots,n$.
\end{lkad_lemma}
\begin{proof}
We prove the statement about a subalgebra having left depth two, namely,  $A \lkado_B A  \lkad| qA$ as $B$-$A$-bimodules. To this apply the additive functor $- \lkado_A {}_{\rho_i}\lkadK$, which results
in $A/ AB_i^+ \lkad| q\lkadK$ as left $B$-modules.  The annihilator of $q\lkadK$ restricted to $B$ is of course
$B^+_i$, which then also annihilates $A/ AB_i^+ $, so $  B_i^+A\subseteq AB_i^+$.   This holds
for each $i = 1,\ldots,n$.
The opposite inclusion is similarly shown to be satisfied by a right depth $2$ extension of augmented algebras.
\end{proof}
 The next theorem computes the depth $d(B,A)$ of the top subalgebra $A / \lkadrad A \cong \lkadK^n$,
or subalgebra of diagonal matrices,
in the path algebra $A$ of an acyclic quiver as given in (\ref{lkadison/eq: kaycue}).
\begin{lkad_theorem}
Suppose the number of vertices $n >1$ in the quiver $Q$, $A = \lkadK Q$ and $B = {\lkadK}^n$. Then depth $d(B,A) = 3$.
\end{lkad_theorem}
\begin{proof}
If the subalgebra in question has depth 1, it has depth 2.  But if it has left depth 2, the lemma above
applies, so that $B_i^+ A \subseteq AB_i^+$ for each $i = 1,\ldots,n$. Note that $AB_i^+$ are
all the lower triangular matrices of the form in (\ref{lkadison/eq: kaycue}) having only $0$'s on column $i$;
similarly, $B_i^+A$ are the triangular matrices having only zeroes on row $i$.  It follows that
$\lkadeps_j A \lkadeps_i = 0$ for each $j = i+1,\ldots, n$.  But $\lkadeps_j (\lkadK Q) \lkadeps_i$ consists of
all the paths from $j$ to $i$.  Since this holds for each $i$, $Q$ consists of $n$ points with no edges; thus we have contradicted the assumption
that $Q$ is connected.  The same contradiction is reached assuming
$B \subset A$ has right depth $2$.

Next it is shown that ${}_BA \otimes_B A_B$ divides a multiple of ${}_BA_B$.
Let $\dim \lkadeps_i A \lkadeps_j = n_{ij}$ .  Then it is clear from (\ref{lkadison/eq: kaycue}) and simple matrix arithmetic that ${}_BA_B \cong \oplus_{n\geq i \geq j \geq 1} n_{ij} \lkadK_{ij}$.

Now $$A \otimes_B A =
\oplus_{i,j = 1}^n \! \oplus_{i \geq k \geq j} \, \lkadeps_i A \lkadeps_k \otimes_B \lkadeps_k A \lkadeps_j $$
since each $\lkadeps_j \in B$ and for each $r \neq k$, $\lkadeps_k \lkadeps_r = 0$.  It follows
that ${}_BA \otimes_B A_B \cong \oplus_{n \geq i \geq j \geq 1} m_{ij} \lkadK_{ij}$ where
$m_{ij} = \sum_{i \geq k \geq j} n_{ik}n_{kj}$.  Since $n_{ii} = 1$ for each $i$, it follows that
$m_{ij} \geq n_{ij}$; moreover, $n_{ij} = 0$ implies $m_{ij} = 0$, since otherwise there is
a path from $i$ to $j$ via some $k$ such that $i \geq k \geq j$.

From the last remark it follows that there is  $q \in \lkadN$ such that $A\otimes_B A \lkad| qA$
as $B$-$B$-bimodules. Thus the minimum depth $d(B,A) = 3$.
\end{proof}
\section{Depth of arrow subalgebra in acyclic quiver algebra}

In this section we compute the depth of the primary arrow subalgebra $B = \lkadK 1_A \oplus A_1 \oplus A_2 \oplus \cdots = \lkadK 1_A + \lkadrad A$ in the path algebra $A$ of an acyclic quiver $Q$, which is of the form
\begin{equation}
\label{lkadison/eq: kaycued}
A = \left( \begin{array}{cccc}
\lkadK & 0 & \cdots & 0 \\
\lkadeps_2(\lkadK Q)\lkadeps_1 & \lkadK & \cdots & 0 \\
\vdots &  \vdots  &  & \vdots \\
\lkadeps_n (\lkadK Q) \lkadeps_1 & \lkadeps_n (\lkadK Q) \lkadeps_2 & \cdots & \lkadK
\end{array}
\right)
\end{equation}
Note that $B$ is a local algebra and augmented algebra with one augmentation $\epsilon: B \rightarrow \lkadK$
equal to the canonical quotient map $B \rightarrow B / \lkadrad B \cong \lkadK$.  We denote the $B$-simple
by $\lkadK_{\epsilon}$ as a pullback module. Again there are $n$ augmentations of $A$ denoted by
$\rho_i$ defining $n$ simple $A$-$B$-bimodules denoted by ${}_i\lkadK_{\epsilon}$, $i = 1,\ldots,n$.

\begin{lkad_lemma}
The natural $B$-$B$-bimodule $A$ is indecomposable.
\end{lkad_lemma}
\begin{proof}
It suffices to show that $\lkadEnd {}_BA_B$ is a local ring \cite{ASS,P}.  Let $F \in \lkadEnd {}_BA_B$ and
choose an ordered basis of $A$ given by $I = \lkadbra \lkadeps_1,\ldots,\lkadeps_n,\alpha_1,\ldots, \alpha_m \lkadket$ where the length of the path $\alpha_i$ is less than or equal to the length of $\alpha_{i+1}$,
all $i = 1,\ldots,m-1$.   Consider the matrix  with $\lkadK$-coefficients, $M = (M^{\alpha}_{\beta})_{\alpha, \beta \in I}$ of $F$ relative to $I$;
then $F(\alpha) = \sum_{\beta \in I} M^{\alpha}_{\beta} \beta$.

Given a path of length $r \geq 1$,  $(i | \alpha | j) \in A_r$, note that $F(\alpha) = \alpha F(\lkadeps_j) = F(\lkadeps_i)\alpha$, so that
$$ \sum_{\beta \in I} M^{\alpha}_{\beta} \beta = \sum_{\gamma \in I} M^{\lkadeps_j}_{\gamma} \alpha \gamma = \sum_{\delta \in I} M^{\lkadeps_i}_{\delta} \delta \alpha. $$
It follows that $M^{\lkadeps_j}_{\gamma} = 0$ for paths $(j|\gamma|k)$ and
$M^{\lkadeps_i}_{\delta}=0$ for all paths $(\ell | \delta | i)$.  Also $M^{\alpha}_{\beta} = 0$ for
all path $\beta \not \in \lkadeps_i A \lkadeps_j$, i.e. not a path from $i$ to $j$.  Finally  deduce
that $M^{\alpha}_{\beta} = 0$ if $\beta \in \lkadeps_i A \lkadeps_j$ but $\beta \neq \alpha$
and $M^{\alpha}_{\alpha} = M^{\lkadeps_i}_{\lkadeps_i} = M^{\lkadeps_j}_{\lkadeps_j}$.

For $i \neq j$ and $\alpha \in \lkadeps_k A \lkadeps_i$, note that $\alpha F(\lkadeps_j) = F(\alpha \lkadeps_j) = 0$,
so that $ \sum_{\beta \in I} M^{\lkadeps_j}_{\beta} \alpha \beta = 0$ implies $M^{\lkadeps_j}_{\beta} = 0$
whenever $s(\beta)=i$.  In particular, $M^{\lkadeps_j}_{\lkadeps_i} = 0$.
It follows that the set of $F \in \lkadEnd {}_BA_B$ has the form of a triangular matrix algebra with constant diagonal, like $B$, and is a local algebra.
\end{proof}

\begin{lkad_theorem}
The depth of the primary arrow subalgebra $B$ in the path algebra $A$ defined above is $d(B,A) = 4$.
\end{lkad_theorem}
\begin{proof}
We first compute $A \otimes_B A$ and show $d(B,A) > 3$.  Note that two paths of nonzero length,
$\alpha$, $\beta$ where $s(\alpha) = i$ satisfy $\alpha \otimes_B \beta = \lkadeps_i \otimes_B \alpha \beta$, which is
zero unless $t(\alpha) = s(\beta)$.  It follows that
$$ A \otimes_B A = \oplus_{i=1}^n \lkadK \lkadeps_i \otimes_B \lkadeps_i \oplus_{i=2}^n \oplus_{j=1}^{i-1} \lkadeps_i \otimes_B \lkadeps_i A \lkadeps_j \oplus_{i\neq j} \lkadK \lkadeps_i \otimes_B \lkadeps_j .$$
It is obvious that the first two summations above are isomorphic as $B$-$B$-bimodules to ${}_BA_B$.
Note that when $i \neq j$, for all paths $\alpha$, $\beta$,
$$\alpha \lkadeps_i \otimes_B \lkadeps_j = 0 = \lkadeps_i \otimes_B \lkadeps_j \beta$$
since $\alpha \lkadeps_i \in B$ is either zero or a path ending at $i$, whence $\alpha \lkadeps_i \lkadeps_j = 0$.
It follows that $A \otimes_B A \cong A \oplus n(n-1){}_{\epsilon}{\lkadK}_{\epsilon}$ as $B$-$B$-bimodules;
moreover, as $A$-$B$-bimodules, we note for later reference
\begin{equation}
\label{lkadison/eq: tensor-square}
{}_A A \otimes_B A_B \cong {}_AA_B \oplus \oplus_{i=1}^n (n-1) {}_i{\lkadK}_{\epsilon}
\end{equation}
By lemma, ${}_BA_B$ is an indecomposable, but the $B$-$B$-bimodule $A \otimes_B A$ contains
another nonisomorphic indecomposable, in fact ${}_{\epsilon}{\lkadK}_{\epsilon}$, so that as $B$-bimodules, $A \otimes_B A \oplus * \not \cong qA$ for any multiple $q$ by Krull-Schmidt.

Now we establish that the subalgebra $B\subseteq A$ has right depth $4$ by comparing (\ref{lkadison/eq: tensor-square}) with the computation below:
$$ A\otimes_B A \otimes_B A = \oplus_{i=1}^n \ \lkadK \lkadeps_i \otimes \lkadeps_i \otimes \lkadeps_i \oplus_{i=2}^n \! \oplus_{j=1}^{i-1} \  \lkadeps_i \otimes \lkadeps_i \otimes \lkadeps_i A \lkadeps_j \oplus_{i \neq j \neq k}\ \lkadK \lkadeps_i \otimes \lkadeps_j \otimes \lkadeps_k $$
$$ \cong A \oplus (n^2 - 1)\, {}_1{\lkadK}_{\epsilon} \oplus \cdots \oplus (n^2 -1)\, {}_n{\lkadK}_{\epsilon}$$
as $A$-$B$-bimodules, where $i \neq j \neq k$ symbolizes $i \neq j$, $j \neq k$ or $i \neq k$.  It is clear that since no new bimodules appear in a decomposition of
${}_A A \otimes_B A \otimes_B A_B$ as compared with ${}_AA \otimes_B A_B$, that there
is $q \in \lkadN$ (in fact $q = n+1$ will do) such that $A \otimes_B A \otimes_B A \lkad| qA \otimes_B A$
as $A$-$B$-bimodules.  It follows that the minimum depth $d(B,A) = 4$.
\end{proof}
It is easy to see from the proof that as natural $B$-$A$ bimodules $A \otimes_B A \otimes_B A \lkad| (n+1)A \otimes_B A$
for very similar reasons.  Note the general fact that  ${}_AA_B$ or ${}_BA_A$ are indecomposable modules if
$\lkadEnd {}_AA_B \cong A^B$, the centralizer subalgebra of $B$ in $A$, is a local algebra.

\section{Concluding Remarks}

It is well-known and easily computed from (\ref{lkadison/eq: kaycued}) that the path algebra $\lkadK Q$ of the quiver $$Q:\ \  n \longrightarrow n-1 \longrightarrow \cdots \longrightarrow 2 \longrightarrow 1$$
is the lower triangular matrix algebra $T_n(\lkadK)$.  Then we have shown above that
for the subalgebras $B_1 = D_n(\lkadK)$ equal to the set of diagonal matrices, and $B_2 = U_n(\lkadK)$
defined by
\begin{equation}
U_n(\lkadK) = \{ \left( \begin{array}{ccccc}
a & 0 & 0 & \cdots & 0 \\
a_{21} & a & 0 & \cdots & 0 \\
a_{31} & a_{32} & a & \cdots & 0 \\
\vdots & &       &      &            \vdots \\
a_{n1} & a_{n2} & a_{n3} & \cdots & a
\end{array}
\right) | a, a_{ij} \in \lkadK \}
\end{equation}
the depths are given by $d(D_n(\lkadK), T_n(\lkadK) ) = 3$ and $d(U_n(\lkadK), T_n(\lkadK)) = 4$.  Both
are not dependent on the order $n$ of matrices.

This situation is different for another interesting series of subalgebras within $T_n(\lkadK)$ given by
\begin{equation}
J_n(\lkadK) = \{ \left( \begin{array}{ccccc}
a_1& 0 & 0 & \cdots & 0 \\
a_2 & a_1 & 0 & \cdots & 0 \\
a_3 & a_2 & a_1 & \cdots & 0 \\
\vdots & &       &      &            \vdots \\
a_n & a_{n-1} & a_{n-2} & \cdots & a_1
\end{array}
\right) | a_1,\ldots,a_n \in \lkadK \}
\end{equation}
also known as the Jordan algebra. This is isomorphic as algebras to $\lkadK[x]/ (x^n)$, a Gorenstein
dimension zero local ring.  Notice that $U_2(\lkadK) = J_2(\lkadK)$, so $$d(J_2(\lkadK), T_2(\lkadK)) = 4.$$

The interesting fact worth mentioning here is that $d(J_3(\lkadK), T_3(\lkadK)) \geq 6$.  This is based
on computations comparing $A \otimes_B A$ and $A \otimes_B A \otimes_B A$ as $B$-$B$-bimodules,
since a new $2$-dimensional indecomposable turns up in the tensor-cube of the ring extension.

The following seems to be an interesting problem not accessible by the techniques of the previous sections:
\begin{equation}
d(J_n(\lkadK), T_n(\lkadK)) = \ ?
\end{equation}

\subsection{Acknowledgements}  The authors would like to thank Sebastian Burciu for visiting Porto in May 2012 and discussing topics related to this paper.  Research on this paper was partially funded by the European Regional Development Fund through the programme {\small COMPETE}
and by the Portuguese Government through the FCT  under the project
 PE-C/MAT/UI0144/2011.

\bibliographystyle{spmpsci}

\end{document}